\documentclass[journal,twoside,web]{ieeecolor}
\usepackage{lcsys}
\usepackage{mathrsfs}
\usepackage{textcomp}
\def\BibTeX{{\rm B\kern-.05em{\sc i\kern-.025em b}\kern-.08em
    T\kern-.1667em\lower.7ex\hbox{E}\kern-.125emX}}
\usepackage{cite}
\usepackage{amsmath,amssymb,amsfonts,mathtools}
\usepackage{amsthm}
\usepackage{algorithmic}
\usepackage{graphicx}
\usepackage{float}
\usepackage{textcomp}
\usepackage{xcolor}
\newtheorem{remark}{Remark}

\newtheorem{theorem}{Theorem}
\newtheorem{definition}{Definition}
\newtheorem{problem}{Problem}
\newtheorem{corollary}{Corollary}
\DeclareMathOperator{\Tr}{Tr}
\DeclareMathOperator*{\argmin}{arg\,min}
\usepackage{bm}
\usepackage{comment}
\begin{document}

\title{Convex Data-Driven Contraction With Riemannian Metrics}

\author{Andreas Oliveira, \IEEEmembership{Graduate Student Member, IEEE},
Jian Zheng, \IEEEmembership{Graduate Student Member, IEEE},\\
and Mario Sznaier, \IEEEmembership{Fellow, IEEE}
\thanks{Manuscript received 18 February 2025; revised 28 April 2025; accepted 19 May 2025. \textit{(Corresponding author: Andreas Oliveira).}}
\thanks{This work was partially supported by NSF grants and CNS--2038493 and CMMI 2208182, AFOSR grant FA9550-19-1-0005, ONR grant N00014-21-1-2431  and DHS grant 22STESE0001-02-00. The authors  are with the Robust Systems Lab,  ECE Department, Northeastern University, Boston, MA 02115. (e-mails: \{franciscodemelooli.a, \;zheng.jian1\}@northeastern.edu, msznaier@coe.neu.edu)}
}

\pagestyle{empty}
\maketitle
\thispagestyle{empty}

\begin{abstract}
The growing complexity of dynamical systems and advances in data collection necessitate robust data-driven control strategies without explicit system identification and robust synthesis. Data-driven stability has been explored in linear and nonlinear systems, often by turning the problem into a linear or positive semidefinite program. This paper focuses on contractivity, which refers to the exponential convergence of all system trajectories toward each other under a specified metric. Data-driven closed-loop contractivity has been studied for the case of weighted $\ell_2$-norms and assuming nonlinearities are Lipschitz bounded in subsets of $\mathbb{R}^n$. We extend the analysis by considering Riemannian metrics for polynomial dynamics. The key to our derivation is to leverage the convex criteria for closed-loop contraction and duality results to efficiently check infinite dimensional membership constraints. Numerical examples demonstrate the effectiveness of the proposed method for both linear and nonlinear systems.
\end{abstract}


\section{Introduction}

\IEEEPARstart{D}{ata-driven} control (DDC) methods seek to avoid the computational complexity and potential conservatism entailed in traditional plant identification-robust control design pipeline by  synthesizing controllers directly from experimental data.  For Linear Time-Invariant (LTI) systems a number of computationally tractable  DDC techniques have been developed in the past few years  \cite{Allgower-2,de2019formulas,de2021low}. DDC of nonlinear systems is considerably more challenging. Local stabilization can be achieved by estimating the Jacobian of the system around a point combined with LTI data-driven techniques \cite{LL}. Data-driven global  stabilization of  nonlinear systems is a much harder problem. Without assuming some form of nonlinearity, the problem of obtaining stability certificates becomes intractable. To address these difficulties, \cite{dai2020semi} assumed polynomial dynamics and resorted to finding a data-driven density function so that the resulting problem was transformed into linear programs. Similarly, in \cite{Persis-Nonlinearn-DDC}, a positive semidefinite program is derived by using quadratic-like Lyapunov functions to ensure closed-loop stability.

This paper focuses on deriving data-driven control strategies for contractivity under Riemannian metrics. Introduced by the seminal paper \cite{LOHMILLER1998683}, contractivity can be seen in some ways as a ``stronger'' property than stability. For excellent treatments of contractivity with $\ell_{p}$ norms see \cite{Sontag-Aminzare}, \cite{FB-CTDS}. Informally, contractive systems are those where, a suitably defined, distance between different trajectories tends to zero exponentially fast, regardless of the initial conditions.
Contractive systems are attractive for many reasons. Firstly, if a time-invariant system is exponentially contracting in $\mathbb{R}^n$ it implies existence of a global exponentially stable equilibrium point \cite{Aylward_2008}. Secondly, tools from semialgebraic geometry enable the search over a wide class of contraction metrics via semidefinite programs \cite{Aylward_2008}. Thirdly contraction guarantees entrainment to periodic forcing \cite{LOHMILLER1998683,FB-CTDS}. Finally, for constant Riemannian Metrics contraction is invariant under any additive disturbance \cite{AlonInvarianceAdditiveDisturbance}.

A convex optimization approach to finding control contraction metrics in a model-based setting was first studied in \cite{manchester2017control}. The approach developed there is the basis of our data-driven results. To enhance its computational efficiency, several extensions were proposed, most notably those utilizing neural network-based contraction metrics. For the case of known dynamics \cite{sun2021learning,Tsukamoto_2021} propose the use of a Neural Network (NN) to co-synthesize a contraction metric and a controller. Unknown dynamics were addressed in 
\cite{chou2021model}, where a NN is used to first identify a model and then learn a   contraction metric $M$. This is followed by a validation step to find a trusted domain $\mathcal{D}$ along with bounds on the identification error, and to verify that $M$ is indeed a contraction in $\mathcal{D}$. While this approach was successfully used to control several robotic systems, it suffers from the conservatism inherent in the model error bounding step and the non-convexity of the two step process.

Alternatively, \cite{yi2023equivalence} exploits the equivalence between Koopman and contraction approaches for systems with a hyperbolic equilibrium point to learn contraction metrics from data for non-actuated systems. However,  dealing with  actuated systems generically requires bilinear, rather than linear, immersions \cite{Otto2024}. Further, if the open-loop plant has  multiple omega-limit sets, continuous one-to-one linear immersions do not exist \cite{2024_koopman_journal}. Finally, \cite{hu2024enforcingcontractiondata} is the closest to the present paper, in the sense that it finds a constant matrix $P$ such that
$\delta x^{\top}P^{-1}\delta x$ is a  control contraction metric for all systems compatible with the observed data and priors. However, as illustrated  in the examples in Section V, limiting the search to state-independent metrics can lead to infeasible problems, even if a state-dependent contraction metric exists.

\textit{Contributions: }To address the issues noted above, we introduce a
 novel convex data-driven formulation for finding Riemannian control contraction metrics 
directly from noisy data. Specifically, our contributions are:
\begin{enumerate}
    \item  A framework for non-conservatively recasting the data-driven control contraction problem, where the Riemannian metrics are restricted to the second-order polynomial matrices, into a sum-of-squares (SOS) problem. By working directly with an exact characterization of the set of all plants consistent with the observed data and priors, this formulation avoids  the conservatism entailed in overbounding this set or the noise description.
    \item  A reformulation of the above problem as a functional Linear Program (LP), which via duality is independent of the parameters of the unknown plant.
    \item Recasting this functional LP as a finite-dimensional Semi-Definite Program (SDP). Notably, the maximal size of the matrices in this SDP is at least $\mathcal{O}(\text{(order of the system)}^2)$ smaller than those appearing in the original SOS problem, resulting in substantial computational complexity reduction.
\end{enumerate}

The paper is organized as follows: Section II provides the definitions and theorems that will be used throughout the paper, along with the problem definition. Section III presents a sufficient condition that solves the problem exploiting Lagrange duality. Section IV develops a tractable relaxation of the dual optimization problem. Section V illustrates these results  with numerical examples. Finally,  Section VI presents  conclusions and directions for future research.

\section{Preliminaries}

\subsection{Notation}

\vspace{-2em}

\begin{alignat*}{2}
    & \mathbb{R}^n (\mathbb{R}_+^n) \quad&&  n\text{-tuples of (non-negative) real numbers} \\
    & \bm{1}, \bm{0}, \bm{I} \quad && \text{vector/matrix of all 1s, 0s, identity matrix} \\
    & \|x\|_{p} && \ell_{p} \text{-norm of vector $x$} \\
    & X^{\top} && \text{the transpose of } X \\
    & x \succeq 0 && \text{$x$ is element-wise non-negative} \\
    & (x)_k && \text{denotes the } k_{\text{th}} \text{element of a vector} \\
    & S_n^{++} && \text{the set of positive definite matrices} \\
    & S_n^{++}[x] && \text{positive definite matrix-valued polynomials} \\
    & X \succeq 0  && X \text{ is positive semidefinite} \\
    & X \succ 0 && X \text{ is positive definite} \\
    & f \in C^d && \text{the $d^{th}$ derivative of $f$ exists and is continuous} \\
    & \text{vec}(X) && \text{vectorization of matrix $X$ along columns} \\
    & \otimes && \text{matrix Kronecker product}  \\
    & \Tr(A) && \text{trace of a matrix $A$}
\end{alignat*}

\subsection{Kronecker Product Property}

The following property of Kronecker product \cite{Matrix-Analysis} will be used in the paper:
\begin{equation}\label{eq:kronecker-property}
    \text{vec}(B^{\top} X^{\top} A^{\top}) = (A \otimes B^{\top})\text{vec}(X^{\top}).
\end{equation}

\subsection{Sum-of-Squares}\label{sec:SOS}
Here we briefly review some concepts from semialgebraic optimization that will be used to recast the data-driven contraction problem into a tractable convex optimization. A polynomial $p(x)$ is an SOS iff it can be written as $p(x)=\sum_i q_i(x)^2$, for some polynomials $q_i$. Alternatively, $p(x)$ is SOS iff there exists a polynomial vector $v(x)$ and a matrix $Q\succeq 0$ such that
$p(x)=v(x)^{\top}Qv(x)$.

A semialgebraic set $\mathbb{K}$ is defined by a finite number of bounded degree polynomials $\{g_i(x)\}_{i=1}^{N_g}$ and $\{h_j(x)\}_{j=1}^{N_h}$:
\begin{equation}\label{eq:bsa}
\begin{aligned}
& \mathbb{K} = \{x\in \mathbb{R}^n \mid  g_i(x) \geq 0, \  h_j(x) = 0, \\ 
& \qquad  \qquad i = 1,\dots,N_g, \; j = 1,\dots, N_h\}.
\end{aligned}
\end{equation}
The set $\mathbb{K}$  {is \textit{Archimedean}} if there exists an $R > 0$ such that {the polynomial $R- \|x\|_2^2$} satisfies  \[R - \|x\|_2^2 =  \sigma_0(x) + \textstyle \sum_i {\sigma_i(x)g_i(x)} + \textstyle \sum_j {\phi_j(x) h_j(x)}\] for some polynomials $\phi_j(x)$ and SOS polynomials $\sigma_i(x)$.
A matrix-valued polynomial $P \in S^{++}_n[x]$ if $\forall x \in \mathbb{R}^n: P(x) \in S^{++}_n$. A matrix-valued polynomial is SOS if there exists a vector of polynomials $v(x)$ and a \textit{Gram} matrix $Q \succeq  0$ such that  $P(x) = (v(x) \otimes \bm{I})^{\top} Q (v(x) \otimes \bm{I})$. 
Every positive definite matrix-valued
polynomial $P(x)$ over an Archimedean set $\mathbb{K}$ satisfies (Theorem 2 of \cite{scherer2006matrix}, Scherer Positivstellensatz)
\begin{equation}
P(x) = \sigma_0(x) + \textstyle \sum_i {\sigma_i(x)g_i(x)} + \textstyle \sum_j {\phi_j(x) h_j(x)}. \label{eq:psatz_noeps}    
\end{equation}
for 
$\sigma_i(x)$ SOS polynomial matrices and
$\phi_j$ symmetric polynomial matrices.
The Putinar Positivstellensatz is the restriction of the Scherer Positivstellensatz to $n=1$ \cite{putinar1993compact}.

\subsection{Contraction Theory}
Informally, a system described by an ODE of the form
\begin{equation}\label{eq:autonomous-ode}
    \dot{x}(t) = f(x), \quad f \in C^1(\mathbb{R}^n,\mathbb{R}^n)
\end{equation}
is contractive if there exists a uniformly positive definite metric tensor, meaning $M(x) \in C^1(\mathbb{R}^n,S_n^{++})$ and $\alpha \bm{I} \preceq M(x)$ for some $\alpha > 0$, such that  the associated Riemannian
differential distance between trajectories decreases with time, that is  $\frac{d (\delta x^{\top} M(x)\delta x)}{dt}<0$. This concept leads to the following definition:

\begin{definition}[\hspace{1sp}\cite{Aylward_2008}]\label{def:contract}
The system \eqref{eq:autonomous-ode} is said to be exponentially contracting, if there exists a uniformly positive definite metric $M \in C^1(\mathbb{R}^n,S_n^{++})$, a  real number $\lambda > 0$ such that the following matrix inequality holds $\forall x \in \mathbb{R}^n$:
\begin{equation} \label{eq:contraction}
    \begin{aligned}
        \frac{\partial f(x)}{\partial x}^{\top}M(x) + M(x)\frac{\partial f(x)}{\partial x} + \dot{M}(x) \prec -2\lambda M(x).
    \end{aligned}
\end{equation}
\end{definition}
In the sequel, by a slight abuse of notation we will refer to systems satisfying \eqref{eq:contraction} as contractive. Theorem 1 of \cite{Aylward_2008} shows  that if \eqref{eq:contraction} is satisfied then there exists a single global exponentially stable equilibrium point.

\subsection{The Data-Driven Contraction Problem}

In this paper, we will consider continuous-time control-affine polynomial systems of the form
\begin{equation}
\label{eq:dynamics-non-autonomous}
    \dot{x} = f(x) + Gu =  F\phi(x) + Gu 
\end{equation}
where $x \in \mathbb{R}^n$ and $u \in \mathbb{R}^m$ are state and control, $f$ is a polynomial up to certain degree, $\phi(x)$ represents a vector of monomials of $x$ and $F, G$ are constant matrices.  For example, 
\begin{equation}
    f(x) = \begin{bmatrix}
        3x_2 - x_1^2 \\
        x_2
    \end{bmatrix} = 
    \begin{bmatrix}
        3 & -1 \\
        1 & 0
    \end{bmatrix}
     \begin{bmatrix}
         x_2 \\
         x_1^2
     \end{bmatrix} = F \phi(x).
\end{equation}
Assume that $T$ noisy measurements $\left\{ \dot{x}[i],x[i],u[i] \right\}_{i=1}^{T}$ satisfying:
\begin{equation}\label{eq:measurements}
 \dot{x}[i] = f(x[i]) + Gu[i] + \eta[i], \|\eta[i]\|_\infty \leq \epsilon,  \forall i = 1,...,T
\end{equation}
 are available. Here the $\ell_\infty$ bounded noise $\eta$ models, for instance, process disturbances or the error incurred when approximating $\dot{x}$ by finite differences.
\begin{definition} The consistency set $\mathcal{P}_1$  of the system \eqref{eq:dynamics-non-autonomous} is the set of all $F$ compatible with the measurements \eqref{eq:measurements} in the following sense:
\begin{equation}\label{eq:MariosFormulation1}
    \begin{aligned}
        \mathcal{P}_1 \doteq  \{F\colon \| F\phi(x[i]) + Gu[i] -\dot{x}[i] \|_{\infty}  \leq \epsilon,\; i = 1,...,T\} 
    \end{aligned}  
\end{equation}
\end{definition}

\begin{problem} \label{prob:problem1} Given noisy data $\left\{ \dot{x},x,u \right\}$ generated by a system of form \eqref{eq:dynamics-non-autonomous}, with a known $G$, find a smooth $M(x)$ and a state feedback control law $u(x)$ such that for all  $F \in \mathcal{P}_1$, the closed-loop system is contractive under the metric induced by $M(x)$.
\end{problem}
\begin{remark} The case where $G$ is a potentially unknown polynomial function of $x$  can be reduced to the problem above by filtering the control action with known dynamics and absorbing the unknown $G(x)$ into $f(x)$, e.g.,
\begin{equation}\begin{aligned}
 \begin{bmatrix}\dot{x} \\ \dot{u} \end{bmatrix} = \begin{bmatrix}f(x) +G(x)u \\
 \bm{0} \end{bmatrix} +  \begin{bmatrix}\bm{0} \\ \bm{I}  \end{bmatrix}v.
\end{aligned}
\end{equation}
\end{remark}

\section{Data-Driven Contraction Control}\label{sec:DD}

The goal of this section is to establish tractable conditions for finding a (differential) data-driven control law that renders all systems in $\mathcal{P}_1$ contractive. This will be accomplished by recasting the problem into a robust optimization problem.

\subsection{Robust Optimization Reformulation}\label{sec:reformulation}
Begin by rewriting the consistency set as 
\begin{equation}\label{eq:MariosFormulation}
    \begin{aligned}
        & \mathcal{P}_1 \doteq \{F\colon \text{Tr}(F\Phi_{i,k}^{\pm}) \leq d_{i,k}^{\pm}, \\
        & \qquad  \qquad \forall i=\{1,\ldots,T\},k=\{1,\ldots,n\} \}
    \end{aligned}  
\end{equation}
where $\Phi_{i,k}^{\pm}$ is a matrix with $\pm \phi(x[i])$ in its $k_{\text{th}}$ column and zeros elsewhere:
\begin{equation}
    \begin{aligned}
    \Phi_{i,k}^{\pm} = 
    \begin{bmatrix}
        \bm{0}, ..., \pm \phi(x[i]), ..., \bm{0}
    \end{bmatrix}, \\
     \text{ and } d_{i,k}^{\pm} = \epsilon \pm (\dot{x}[i] - Gu[i])_{k}.
   \end{aligned}
\end{equation}

For example, for a single measurement $x[1],u[1] \in$ $\mathbb{R}^2$, we have
\begin{equation}
\begin{array}{ll}
    \Phi_{1,1}^{\pm} = 
    \begin{bmatrix}
        \pm\phi(x[1]) & \bm{0}
    \end{bmatrix},
    & d_{1,1}^{\pm} = \epsilon \pm (\dot{x}[1]-Gu[1])_1 \\
    \Phi_{1,2}^{\pm} = 
    \begin{bmatrix}
         \bm{0} & \pm \phi(x[1])
    \end{bmatrix},
    & d_{1,2}^{\pm} = \epsilon \pm (\dot{x}[1]-Gu[1])_2.
\end{array}
\end{equation}
In the sequel the $\pm$ superscript in $\Phi_{i,k}^{\pm}$ and $d_{i,k}^{\pm}$ is omitted to keep the notation cleaner.

For actuated systems of the form
\eqref{eq:dynamics-non-autonomous}, condition \eqref{eq:contraction} now depends on finding a differential control input $\frac{\partial u(x)}{\partial x}$ that satisfies:
\begin{equation}\label{eq:contraction1}
\begin{aligned}
  & \dot{M}(x) +  \left(F\frac{\partial \phi(x)}{\partial x}   + G\frac{\partial u(x)}{\partial x}\right)^{\top}M(x) \\ & \qquad + M(x)\left(F\frac{\partial \phi(x)}{\partial x} + G\frac{\partial u(x)}{\partial x}\right)  \prec -2\lambda M(x).
\end{aligned}
\end{equation}

If $M(x),\frac{\partial u(x)}{\partial x}$ exist and satisfy the above then $M(x)$ is referred to as a control contraction metric. Using the change of variable $W=M^{-1}$ and leveraging Proposition 2 in \cite{manchester2017control}, the differential control law $\frac{\partial u(x)}{\partial x} = -\frac{1}{2}\rho(x)G^{\top}W^{-1}(x)$ renders the closed-loop system contractive if there exists $W(x) \in S_n^{++}$ and a function $\rho(x)$ such that $\forall x \in \mathbb{R}^n$ (we omit $(x)$ for brevity hereafter):

\begin{equation}\label{eq:contraction-closed-loop}
    \begin{aligned}
         -\dot{W} + W(F\frac{\partial \phi}{\partial x})^{\top} + (F\frac{\partial \phi}{\partial x})W  + 2\lambda W
         -\rho GG^{\top} \prec 0
    \end{aligned}
\end{equation}

\begin{equation}\label{eq:G-orthogonal-to-metric}
\sum_j   \frac{\partial W(x)}{\partial x_j} (Ge_i)_j = 0,\; \forall i = 1,...,n
\end{equation}
where $e_i$ is the $i$-th basis vector in $\mathbb{R}^n$. For the future, condition (\ref{eq:G-orthogonal-to-metric}) is abbreviated as $\partial_G W(x) = 0$.

In terms of \eqref{eq:MariosFormulation}, \eqref{eq:contraction-closed-loop} and \eqref{eq:G-orthogonal-to-metric}, Problem \ref{prob:problem1} can be reformulated as:
\begin{problem}\label{prob:RO}
Find $W \in S^{++}_n[x]$, $\rho$ and $\lambda > 0$ such that \eqref{eq:contraction-closed-loop}-\eqref{eq:G-orthogonal-to-metric} hold for all $F$ satisfying the trace condition 
\eqref{eq:MariosFormulation}.
\end{problem}

\begin{remark}\label{rem:W} Since $M=W^{-1}$, in principle global contractivity requires $W(x)$ to be uniformly bounded above, e.g., $W(x) \preceq \sigma \bm{I}; \forall x \in \mathbb{R}^n$, which clearly cannot be accomplished with a polynomial matrix. However, Lemma 1 in \cite{manchester2017control} shows that this condition can be relaxed to existence of a quadratic bound on the maximum eigenvalue of $W(x)$, that is $\lambda_{\text{max}} (W(x)) \leq \|Ax + B\|_2^2$ for some matrix $A$ and vector $B$. Since this condition is automatically satisfied for $W \in S^{++}_2[x]$, in the sequel we will consider only second order polynomial matrices.
\end{remark}

\subsection{Solution via Duality}

In principle, Problem \ref{prob:RO} can be reduced to sequence of SDPs by restating it as a polynomial optimization in the indeterminates $x$ and $F_{ij}$, the elements of $F$, and exploiting  Scherer's Positivstellensatz \cite{Scherer-PSATZ}. However, this approach is practically limited to relatively small, low-order systems due to the very poor scaling properties of the resulting SDP with respect to the size of $F$. To avoid this difficulty, in this paper, we will pursue a duality-based approach to obtain an equivalent condition that does not involve $F$.

Assume for now that $x$ is fixed and hence $\rho$ and $W$ are constant. Scalarizing  \eqref{eq:contraction-closed-loop} leads to the following equivalent condition:
\begin{equation} \label{eq:primal}
    \begin{aligned}
 &   y^{\top}(  W\partial\phi^{\top} F^{\top} + F\partial\phi W - \dot{W} + 2\lambda W -\rho GG^{\top})y < 0  \\
& \text{ $\forall \|y\|_2=1$ and $F$ that satisfy}    \Tr(F\Phi_{i,k}) \leq d_{i,k}. \\
    \end{aligned}
\end{equation}
This condition can be reduced to an SDP via Putinar's Positivstellensatz in $y$ and $F_{ij}$. However, as before, this leads to problems with poor scaling properties. Rather than pursuing this approach,  we will enforce \eqref{eq:primal} by computing
\begin{equation} \label{eq:primal2}
    \begin{aligned}
 p^*(y,W,\rho) =  \max_{F} y^{\top}(  W\partial\phi^{\top} F^{\top} + F\partial\phi W - \dot{W} \\ + 2\lambda W -\rho GG^{\top})y   \\
 \text{subject to: }
d_{i,k} -\Tr(F\Phi_{i,k}) \geq 0 \\
    \end{aligned}
    \end{equation}
and finding $W(x)$ and $ \rho(x)$ such that $p^*(y,W,\rho) < 0$ for all $\|y\|_2=1$ and such that \eqref{eq:G-orthogonal-to-metric} is satisfied.
\begin{theorem}\label{teo:dual} A metric tensor $W(x) \in C^1(\mathbb{R}^n, S_n^{++})$ and a function $\rho(x)$ solve Problem 2 
if there exist non-negative functions $\mu_{i,k}(x,y)$ such that:
\begin{subequations}\label{eq:dual_problem}
    \begin{align} 
        &y^{\top}(-2\lambda W + \rho GG^{\top})y  -\sum_{i,k}\mu_{i,k}d_{i,k} >  0 \;  \label{subeq:1}\\
         & \text{$\forall \|y\|_2=1$}  \nonumber \\
         & \text{vec}(-2\partial \phi W yy^{\top} + \sum_{i,k} \mu_{i,k}\Phi_{i,k})^{\top}  \label{eq:zero-constraint} \\
         & \qquad \qquad + \sum_{i} \Tr(\frac{\partial W}{\partial x_i} yy^{\top})(e_i^{\top} \otimes \phi^{\top}) = 0 \nonumber\\
        & \partial W_G(x) = 0 \label{eq:G_orthogonal}\\
        & \mu_{i,k}(x,y) \geq 0,\;  W(x) \succ 0.
    \end{align}
\end{subequations}
Moreover, if the consistency set $\mathcal{P}_1$ has a non-empty interior, then the condition is also necessary. 
\end{theorem}
\begin{proof}
 The Lagrangian of \eqref{eq:primal2} for fixed $y,W,\rho$ is given by:
\begin{equation}
    \begin{aligned}
  &  L(F,\mu_{i,k})  = \sum_{i,k}\mu_{i,k}(d_{i,k}-\Tr(F\Phi_{i,k})) + \\ & y^{\top}(  W\partial\phi^{\top} F^{\top} + F\partial\phi W - \dot{W} + 2\lambda W -\rho GG^{\top})y
    =  \\&\Tr\left \{ F\left (2\partial\phi W yy^{\top} - \sum_{i,k}\mu_{i,k}\Phi_{i,k} \right ) \right \}-\\
   & \Tr\left \{ \left(\dot{W} - 2\lambda W +\rho GG^{\top}\right)yy^{\top} \right \} 
    +\sum_{i,k}\mu_{i,k}d_{i,k}.
     \end{aligned}
\end{equation}

Expanding $\dot{W}$ by leveraging \eqref{eq:G_orthogonal} and \eqref{eq:kronecker-property} yields:
\begin{equation}
\begin{aligned}
  \dot{W} = \sum_{i} \frac{\partial W}{\partial x_i} (\phi^{\top}F^{\top})e_i=  \sum_{i} \frac{\partial W}{\partial x_i} (e_i^{\top} \otimes \phi^{\top})\text{vec}(F^{\top}).
\end{aligned}
\end{equation}

Using the former identity and the fact that $\Tr(AB) = \text{vec}(B)^{\top}\text{vec}(A^{\top})$, the Lagrangian becomes:
\begin{equation}
    \begin{aligned}
    L(F,\mu_{i,k}) = \text{vec}(2\partial \phi W yy^{\top} - \sum_{i,k} \mu_{i,k}\Phi_{i,k})^{\top}\text{vec}(F^{\top}) \\ - \sum_{i}\Tr(\frac{\partial W}{\partial x_i} yy^{\top})(e_i^{\top} \otimes \phi^{\top})\text{vec}(F^{\top}) \\ + \Tr((2\lambda W -\rho GG^{\top}) yy^{\top}) + \sum_{i,k}\mu_{i,k}d_{i,k}.
    \end{aligned}
\end{equation}

Notice that  $L(F,\mu_{i,k})$ is affine in $F$, therefore, the dual function: 
\begin{equation}
    \begin{aligned}
        g(\mu_{i,k}) = \sup_{F} L(F,\mu_{i,k})
    \end{aligned}
\end{equation}
is finite only if \eqref{eq:zero-constraint} holds.
 Therefore:
\begin{equation}
    g(\mu_{i,k}) =  
    \begin{dcases}
        {\Tr}\left \{(2\lambda W -\rho GG^{\top})yy^{\top} \right \}  +\sum_{i,k}[\mu d]_{i,k} \\ 
        \qquad \text{if  (\ref{eq:zero-constraint}) holds and}
        \\[2ex]
        \infty \quad \text{ otherwise } 
    \end{dcases}
\end{equation}
Hence, if there exist  non-negative multipliers $\mu_{i,k}(x,y) \geq 0$ satisfying 
\eqref{eq:dual_problem}, $g(\mu_{i,k})<0$. 
From weak duality \cite{boyd2004convex} it follows that $p^*(y,W,\rho) <0$ for all $\|y\|_2 = 1,\rho(x), W(x)\succ 0$. Thus \eqref{eq:contraction-closed-loop}-\eqref{eq:G-orthogonal-to-metric} hold for all $F$ satisfying 
\eqref{eq:MariosFormulation}. Moreover, if the consistency set $\mathcal{P}_1$ has a non-empty interior, then  \eqref{eq:primal2} is a linear program in $F$ which now satisfies Slater's condition at each $y,W(x),\rho(x)$. Hence  strong duality  holds and $p^*(y,W,\rho) = g^*(\mu_{i,k})$. Thus \eqref{eq:dual_problem} is also necessary.
\end{proof}
\begin{corollary}
    The differential feedback control law $\frac{\partial u}{\partial x} = -\frac{1}{2}\rho(x)G^{\top}W^{-1}(x)$ renders all systems in $\mathcal{P}_1$ contractive.
\end{corollary}

\section{Tractable Relaxations}

From Theorem \ref{teo:dual} it follows that Problem \ref{prob:RO} reduces to a feasibility problem in $W(x) \in S_n^{++}, \mu_{i,k}(x,y) \in \mathbb{R}_+$, and $\rho(x)$. However, searching for a matrix function $W(x) \succ 0$ is generically intractable. To avoid this problem we will restrict the search to SOS matrices, that is $W(x) = (\bm{ \Psi \otimes I})^{\top}Q (\bm{\Psi \otimes I})$ where $Q \succ 0$ and  $\bm{\Psi}$ is a basis of monomials per \cite{Scherer-PSATZ}. Further, in view of Remark \ref{rem:W}, we will only consider second order polynomial matrices, that is
 $W(x) \in S_2^{++}[x]$. Similarly, to guarantee that $\mu_{i,k} \geq 0$ the search is performed for $\mu_{i,k}$ belonging to the set of SOS polynomials. Finally, we also restrict $\rho$ to be a polynomial. Under these assumptions, \eqref{eq:dual_problem} is a semi-algebraic optimization over $W,\mu$ and $\rho$, that satisfy the linear constraints in \eqref{eq:zero-constraint},\eqref{eq:G_orthogonal}. Therefore, we utilize Putinar's Positivstellensatz and other standard SOS tools to solve it via a sequence of SDPs.

It is important to get an understanding of the computational complexity before presenting the examples. Assume the system dimension is $n$, the highest degree in $\phi$ is $p$, the highest degree in $W$ is $2q$ (where $q$ is either 0 or 1 in view of Remark \ref{rem:W}) and $T$ data samples were collected. As discussed $W(x) = (\bm{ \Psi \otimes I})^{\top}Q (\bm{\Psi \otimes I})$ where $\bm{\Psi}$ is the vector of all monomials up to degree $q$ of the system state. Hence the size of the Gram matrix $Q$ is  $n{n+q \choose q} \times n{n+q \choose q}$.
 
Furthermore, based on the set $\mathcal{P}_1$ and the number of data samples $T$, there are $2nT$ SOS functions $\mu_{i,k}$ each of which should have degree in $x$ at least $2q+p-1$ and degree 2 in $y$. Thus, the associated Gram matrices have dimension $n_G \times n_G$ with $n_G \approx (n+1){\frac{p-1}{2}+q+n \choose n}$.

On the other hand, enforcing \eqref{eq:contraction-closed-loop} directly through Scherer's Positivstellensatz in the indeterminate $x \in\mathbb{R}^{n}$ and $F\in \mathbb{R}^{n\times l}$, involves an $n \times n$ SOS matrix in $(n+nl)$ variables. Assuming polynomials of order $2q$ in $W(x)$ and $n_F$ in $F$ the corresponding Gram matrix will have dimension $n_S \times n_S$ where $n_S \approx n{\frac{p-1}{2}+q+n\choose n} {nl + n_F \choose n_F}$.
In summary, finding higher-order metrics comes at the cost of combinatorial complexity in the size of the positive semidefinite matrix constraints. However, exploiting duality mitigates the growth of the largest Gram matrix, as compared to straight application of Scherer's Positivstellensatz, when $l>1$, {even when choosing $n_F=1$}. Additionally, while increasing the number of data samples can shrink the consistency set $\mathcal{P}_1$ it increases the amount of SOS functions $\mu_{i,k}$ needed.

\section{Numerical Examples}

To decrease numerical errors and computational complexity, the following examples assume that it is known beforehand which monomials appear in the dynamics. Thus  we will use a dictionary $\phi$ containing only those terms. Additionally due to the limitations of numerical solvers in handling zero equality constraints as required in \eqref{eq:zero-constraint} and \eqref{eq:G_orthogonal}, we adopt the heuristic that any polynomial coefficient smaller than $10^{-5}$ in $W,\rho$ is set to 0. To enforce the other constraints in \eqref{eq:dual_problem} we utilize the \textit{solvesos} function in the YALMIP toolbox and with MOSEK as the underlying solver. Furthermore in the optimization program to solve Problem \ref{prob:RO} we consider $\lambda$ a positive fixed constant. When this choice resulted in infeasibility,  we reran the optimization routine with a smaller $\lambda$ while maintaining the positivity constraint. In practice finding the optimal value for $\lambda$ can be done via a bisection algorithm described in section 2.7 of \cite{FB-CTDS}. Note that a larger $\lambda$ is desirable, since it yields faster convergence to the origin. Each data-driven problem below was solved under 100 seconds, on a Apple M1 Pro with 16GB of RAM. Only 60 measurements were required, namely $\{x[i],\dot{x}[i]\}_{i=1}^{60}$ (the control input $u$ is assumed to be zero), generated with the \textit{ode45} or \textit{ode15s} commands in MATLAB and coming from six randomized initial conditions in the interval $[-1,1]^n$. The derivative $\dot{x}[i]$ is disturbed by adding uniform random noise bounded above by $\|\epsilon\|_{\infty} = \frac{1}{15}\max_i(|\dot{x}[i]|)$. The code is available at: https://github.com/andoliv1/Contraction-For-DDC.
 
\subsection{Linear System}

First, we validate our results in a linear system. Consider the system given by:
\begin{equation}\label{eq:linear1}
    \dot{x} = \begin{bmatrix}
0.4285 & -0.4298 \\
0.4018 & 1.3036
\end{bmatrix}x + \begin{bmatrix}
    -0.7826 & 0.7731 \\
    -0.5110 & 0.0339
    \end{bmatrix}u.
\end{equation}

If $u = [0,0]^{\top}$ the above is unstable with eigenvalues at $0.7291,1.0030$. The objective is to find a controller $u$ and metric tensor $W$ that ensures the system is contractive from data observations. Solving \eqref{eq:dual_problem} produces the following  constant matrix  $W \in S_n^{++}$ and scalar $\rho$

\begin{equation}\label{eq:linear1-contraction}
\begin{aligned}
    W = \begin{bmatrix}
   0.1019 & 0.0154 \\
   0.0154 & 0.0028
    \end{bmatrix},\; \rho =9.9686.
\end{aligned}
\end{equation}

In order to close the loop in the underlying system, one needs to define $u(x)$ from $\frac{\partial u(x)}{\partial x} = -\frac{1}{2}\rho G^{\top}W^{-1}$. Following \cite{leung2017nonlinear}, but simplifying the control input to track the origin:
\begin{equation}\label{eq:integral-control-input}
    u(x) = \int_{0}^{1} -\frac{1}{2}\rho(\gamma(s))G^{\top}W(\gamma(s))^{-1}\frac{\partial\gamma(s)}{\partial s} ds
\end{equation}
where $\gamma: [0,1] \rightarrow \mathbb{R}^2, \gamma(0) = [0,0]^{\top}, \gamma(1) = x$ is the geodesic under the norm induced by $W^{-1}$, which if constant implies $\gamma(s) = sx$. Since both $\rho,W$ are constant we have that the feedback controller that makes the system contractive is $u(x) = -\frac{1}{2}\rho G^{\top}W^{-1}x$. Note that in this case a similar result can be obtained using Krasovskii's criteria \cite{LOHMILLER1998683}.

\subsection{Nonlinear Systems}

We alter the system referenced in \cite{Aylward_2008} such that it is no longer open-loop contractive and it has a control input:
\begin{equation}\label{eq:non-contractive-dynamics-with-input}
    \begin{bmatrix}
        \dot{x}_1 \\
        \dot{x}_2
    \end{bmatrix}
     = 
     \begin{bmatrix}
         -x_2 - \frac{3}{2}x_1^2 - \frac{1}{2}x_1^3 \\
         3x_1 + x_2
     \end{bmatrix}
      + 
      \begin{bmatrix}
          1 & 0 \\
          0 & 1
      \end{bmatrix}u.
\end{equation}

Running the optimization program \eqref{eq:dual_problem} with respect to \eqref{eq:non-contractive-dynamics-with-input} results in:
\begin{equation}
\begin{aligned}
     & W = \begin{bmatrix}
     0.0074  &  -0.0027 \\
    -0.0027  &  0.1018
    \end{bmatrix},
    \\
    & \rho =  3.9508-0.0005x_1+2.9056x_1^2+2.9044x_2^2.
\end{aligned}
\end{equation}

To demonstrate the need for a nonlinear contraction metric, let's draw from the example in \cite{manchester2017control}:
\begin{equation}\label{eq:manchester-ode-3D}
    \begin{bmatrix}
        \dot{x}_1 \\
        \dot{x}_2 \\
        \dot{x}_3
    \end{bmatrix}
     = 
     \begin{bmatrix}
         -x_1 + x_3 \\
         x_1^2 - x_2 - 2x_1x_3 +x_3 \\
         -x_2
     \end{bmatrix}
      + 
      \begin{bmatrix}
          0 \\
          0 \\
          1
      \end{bmatrix}u.
\end{equation}

First it should be noted that trying to find a constant matrix $W$ and polynomial $\rho$ that makes the system closed-loop contractive resulted in infeasible solutions, meaning $W$ was rank deficient or the solver could not find a solution. Therefore, from Remark 2 we should allow $W \in S_2^{++}[x]$. Now running \eqref{eq:dual_problem} the following $W$ and $\rho$ are found that make the closed-loop system contractive (only the (3,3)-th entry of $W$, denoted as $W_{3,3}(x)$ and the first 7 components of $\rho$ are shown for the sake of brevity):
\begin{equation}
\begin{aligned}
    W_{3,3}(x) = 0.0173-0.00027x_1- \\0.0004x_2+0.0015x_1^2+0.0015x_2^2
\end{aligned}
\end{equation}
\begin{equation}
\begin{aligned}
    \rho(x) =2.765 - 0.0002x_1 + 0.3204x_1^2+ \\ 0.2962x_2^2+0.2999x_3^2 \dots
\end{aligned}
\end{equation}

\begin{figure}[H]
    \centering
    \includegraphics[width=0.77\linewidth]{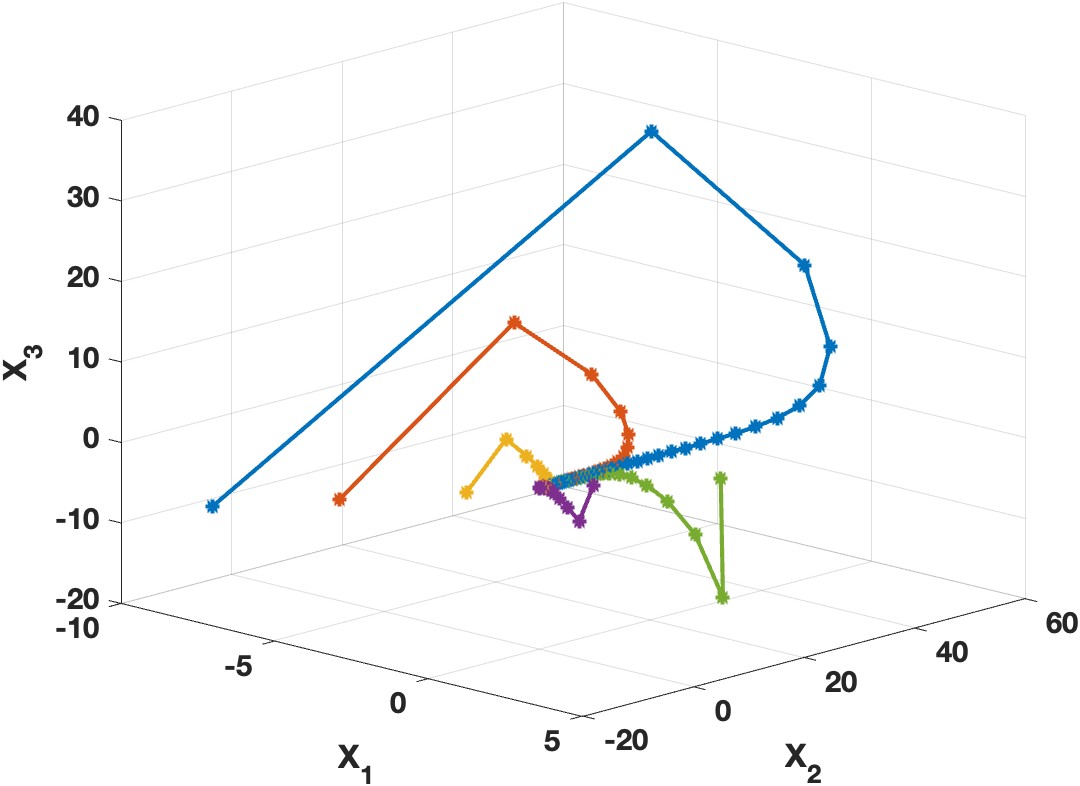}
    \caption{Progression of trajectories for different initial conditions.
    }
    \label{fig1}
\end{figure}

\begin{figure}[H]
    \centering
    \includegraphics[width=0.77\linewidth]{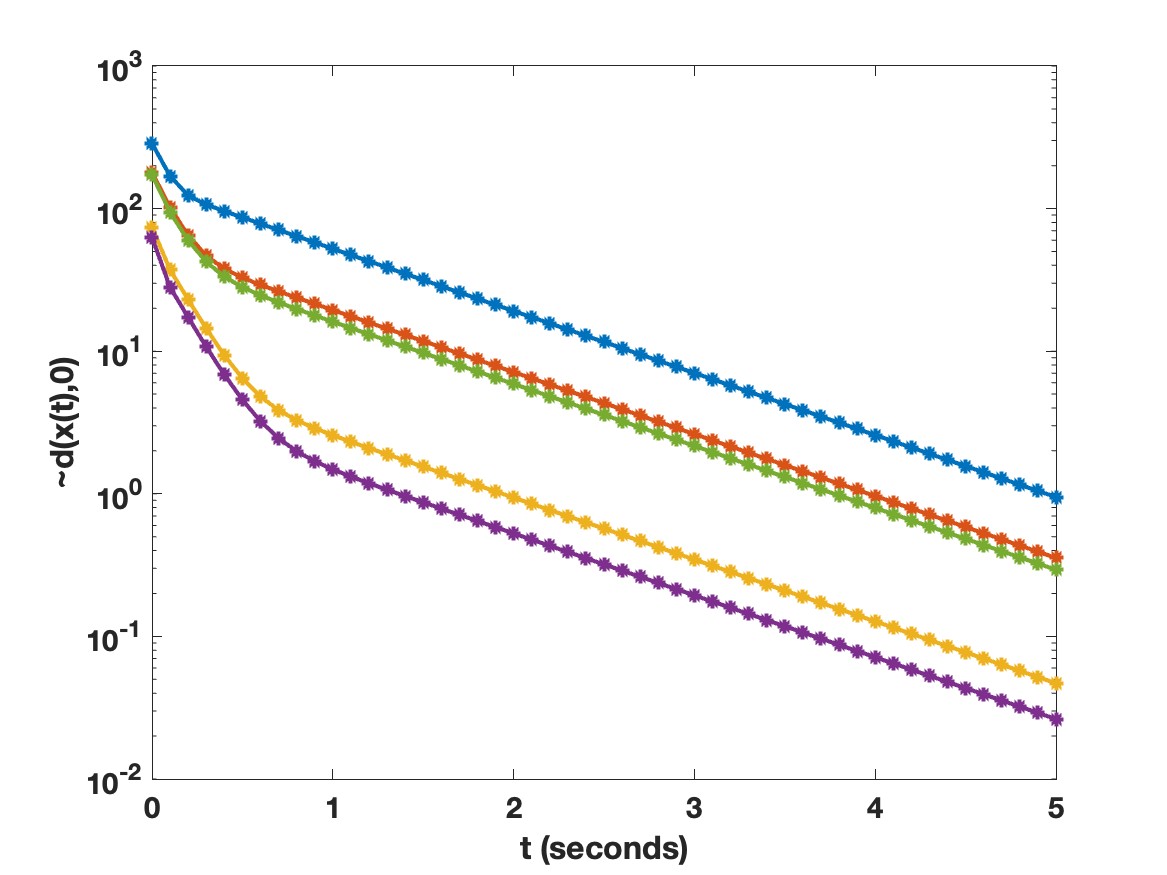}
    \caption{Graph of the approximate distance between $x(t)$ and $\bm{0}$ computed using the approximate geodesic. The piecewise linear relation between $t$ and the y-axis in log scale indicates $x(t)$ is exponentially converging to the origin.}
    \label{fig2}
\end{figure}

As mentioned before, enforcing contraction in the closed loop with respect to the Riemannian metric requires the geodesic in (\ref{eq:integral-control-input}), which amounts to finding:

\begin{equation}\label{eq:geodesic}
    \gamma^* = \argmin_{\gamma \in \Gamma(\bm{0},x(t))} \int_0^1 \frac{\partial \gamma(s)}{\partial s}^{\top}M(\gamma(s))\frac{\partial \gamma(s)}{\partial s} ds
\end{equation}
where $\Gamma(\bm{0},x(t))$ is the set of smooth regular curves that connect the origin to $x(t)$ \cite{manchester2017control}. Finding an explicit solution of the above is computationally hard as noted in \cite{leung2017nonlinear}. Therefore, to obtain an approximate solution we discretize \eqref{eq:geodesic} into a finite sum and turn $\gamma$ into a piecewise linear path with boundary conditions. The order of the discretization of $\gamma$ was set to be $100$, meaning $\gamma$ consisted of $100$ piecewise linear steps. The resulting nonlinear program is fed into \textit{CasADi} using IPOPT as the underlying solver and each run of the problem took less than $0.5$ seconds. Figures \ref{fig1} and \ref{fig2} show the evolution of $5$ different initial value problems of system (\ref{eq:manchester-ode-3D}) with the computed control input applied.

\section{Conclusion}

We study the problem of finding a metric and state feedback controller that renders all systems compatible with noisy experimental data contractive. Previous works on this problem have considered the case of searching for weighted $\ell_2$-norms that render sector-bounded nonlinear systems compatible with data contractive \cite{hu2024enforcingcontractiondata}. Our approach extends this by considering polynomial systems in unbounded domains and searching for Riemannian metrics. The key to our formulation lies in leveraging the convex criteria developed in \cite{manchester2017control} and applying duality to significantly reduce computational complexity. The problem can be solved efficiently for low dimensional systems and low order metric tensor $W$ but remains challenging as the system dimension and order increase. Future work seeks to address this by using alternative characterizations of contractivity based on matrix log norms \cite{FB-CTDS}.

\bibliographystyle{IEEEtran} 
\bibliography{ref}

\end{document}